\documentclass[letterpaper, 10pt, conference]{ieeeconf}
\IEEEoverridecommandlockouts 
\usepackage{amsmath}
\usepackage{amssymb}

\usepackage{amsthm}
\usepackage{amsfonts}
\usepackage{url}
\usepackage{cite}
\usepackage{graphicx}
\graphicspath{{./figures/}}
\newtheorem{theorem}{Theorem}[section]
\newtheorem{lemma}[theorem]{Lemma}

\newcommand{\argmax}{\operatornamewithlimits{arg\,max}}

\title{Multi-Agent Discrete Search with Limited Visibility}
\author{Huanyu Ding$^1$ and David Casta{\~n}{\'o}n$^1$
\thanks{$^1$Department of Electrical and Computer Engineering, Boston University, United States.
        \texttt{\small \{hyding, dac\}@bu.edu}.}
         \thanks{This work was  supported by  NSF award CNS-1330008.}}

\begin{document}
\maketitle

\begin{abstract}
The problem of search by multiple agents to find and localize objects arises in many important applications.  In this paper, we study a class of multi-agent  search problems in which each agent can access only a subset of a discrete search space, with detection performance that depends only on the location.  We show that this problem can be reformulated as a minimum cost network optimization problem, and develop a fast specialized algorithm for the solution.   We prove that our algorithm is correct, and has worst case computation performance that is faster than general minimum cost flow algorithms.  We also address the problem where detection performance depends on both location and agent, which is known to be NP-Hard.  We reduce the problem to a submodular maximization problem over a matroid, and provide an approximate algorithm with guaranteed performance.      We illustrate the performance of our algorithms with simulations of search problems and compare it with other min-cost flow algorithms. 
\end{abstract}

\section{Introduction}

The proliferation of intelligent agents in diverse applications from building security, defense, transportation, and medicine has created a need for automated processing and exploitation of information. An important class of problems related to the applications of these intelligent systems is the detection and localization of  hidden objects of interest, known as search problems \cite{Koopman1946, StoneBook, ahlswede1987search, castanon1995optimal, song2004discrete, ding2015optimal, ding2015multi} in the operation research, statistics and control communities.  Examples of search problems include localizing a lost submarine or plane in the ocean \cite{StoneBook}, determining faulty components in a system \cite{pattipati1990application}, discovering archaeological sites \cite{acc16, ding2016fast}, etc. 

The study of search problems dates back to its early application for objects at sea in the 1940s \cite{Koopman1946}.  in a discrete version of the problem, objects are located in a discrete search space $\mathcal{X}$, which has a finite set of possible locations.  The general framework assumes that there is a single object for the search.  Given a prior probability distribution on the space of locations for where that object is located, a common goal is to select search actions to maximize the probability that the object is detected.

Most of the work in classical discrete search theory \cite{StoneBook, ahlswede1987search, Benkoski1991} considered search with a single agent. An agent's measurements of a location can result in a binary signal indicating whether a detection was obtained or not.  If the object is not present at the location being searched, no detection will be measured (thus, the model does not allow false alarms.)   There are classes of search problems that include more complex sensing models, allowing for false alarms \cite{castanon1995optimal,Davidsbook}.  Those types of search problems result in partially-observed Markov Decision Problems, and often require approximations or special structures to obtain solutions \cite{Davidsbook,markovObj}.  In this paper, we focus on the classical measurement model where missed detections are possible, but not false alarms.

There are usually two types of objectives for the search: detection search and whereabouts search.  In detection search \cite{StoneBook, castanon1995optimal, song2004discrete}, the objective is to maximize the probability of detecting the hidden object subject to the budget on effort.  In whereabouts search \cite{tognetti1968, kadane1971whereabouts, castanon1995optimal}, the objective is to maximize the probability of identifying the correct location of the object subject to the budget on effort.  In this paper, we focus primarily on detection search, although many of the results can be extended to whereabouts search.   


Multi-agent discrete search problems without false alarms were considered by Song and Teneketzis \cite{song2004discrete}.  They assumed that each agent has a budget of $N_{each}$ units of effort and each search action consumes one unit of effort.  In their model, the probability of detection for an agent searching a location was independent of the agent, while depending only on the location.   Every agent could search every location, but there was an additional constraint that no two agents could search the same location at the same time.  The objective was to maximize the probability of finding the object after spending all the budget.  They developed a complete algorithmic solution for constructing an optimal multi-agent search schedule.

In this paper, we generalize the problem in \cite{song2004discrete} to the case where each agent can only search in part of the search space.  We refer to this problem as  {the sparse multi-agent discrete search problem}.  Inspired by \cite{wtaReport}, we propose a novel perspective of viewing the problem as a minimum cost network flow problem.  We propose a fast specialized algorithm for solving the problem based on the min-cost flow perspective.  We show that the algorithm always terminates in finite time, and analyze the time complexity of the algorithm.  We prove that the algorithm yields an optimal agent schedule.  We perform experiments to show that our specialized algorithm is faster than general min-cost flow algorithms such as the capacity scaling algorithm.  

We also consider the general multi-agent sparse search where the probability of detection depends on both the agent performing the search and the location of the search.  This problem has been studied previously and found to be NP-hard  \cite{lloyd1986weapons, ahuja2007wta}.  In \cite{ahuja2007wta}, several approximate algorithms are considered, including branch-and-bound and local search.  We show that this problem can be posed as a submodular maximization problem over a matroid, and establish that the greedy algorithm is guaranteed to perform within a factor of $1/2$ of the optimal value.  We also discuss a more complex algorithm, based on recent ``continuous greedy'' approximations that are guaranteed to perform within a factor of $(1 - 1/e)$ of the optimal value. 

The rest of the paper is organized as follows: In Section \ref{sec:formulation}, we formulate the sparse multi-agent discrete search problem.  In Section \ref{sec:perspective}, we introduce the min-cost flow perspective of viewing the problem.  We provide the primal and dual linear programming formulations of the problem, and derive conditions in which an agent schedule is optimal.  We propose and analyze a fast algorithm for solving the problem in Section \ref{sec:algorithm}, and prove that it yields an optimal solution.  Section \ref{sec:experiments} contains some experiment results.  Section \ref{sec:discussion} discusses the extensions to the case where probability of detection depends both on agent and search location. We conclude the paper in Section \ref{sec:conclusion} with  suggestions for future work.  

\section{Formulation of the sparse multi-agent discrete search problem}
\label{sec:formulation}

Consider $M$ agents that search for a stationary object hidden in one of $K$ discrete locations in discrete time stages.  Each agent has limited search accessibility to only a subset of locations.  Define an \textit{accessibility pair} of agent $m$ and location $k$ if and only if agent $m$ has accessibility to location $k$, denoted by $(m,k)$.  Let $\mathcal{A}$ denote the set of all accessibility pairs.  Without loss of generality, assume that each location can be accessed by at least one agent and each agent can access at least one location.   Assume that a search of any location by any agent that can access it costs one unit of effort.  Let the budget of total effort agent $m$ can allocate be $N_m \in \mathcal{Z}^{+}$, the set of positive integers.  Let $N=\sum_{m=1}^M N_m$ be the total search effort available. 

Denote the prior probability of location $k$ containing the object as $p_{k0}$.  If the object is in location $k$, it will be found by any agent searching the location with probability $\alpha_k$.  If the object is not in location $k$, searching it will always yield ``no detection''.  Assume that agent observations are conditionally independent across locations and of previous searches.  

A useful quantity is $p_{kj}$, the probability that $j-1$ searches of location $k$ have not found the object while the $j$-th search of location $k$ finds the object.  Using the conditional independence of search outcomes, this expression becomes
\begin{align*}
p_{kj} = p_{k0}(1-\alpha_k)^{j-1}\alpha_k
\end{align*}
Our objective is to allocate agent effort among accessible locations so as to maximize the probability of finding the object after spending all budgets.  Suppose that we allocate $x_{mk}$ units of effort to search location $k$ using agent $m$, which sums up to $u_k=\sum_{m=1}^M x_{mk}$ searches for location $k$, so the probability of finding the object after $u_k$ searches of location $k$ is $p_{k0}(1-(1-\alpha_k)^{u_k})$.   Then the problem is to find an agent schedule $\boldsymbol{x} = \{x_{mk}, \forall (m,k) \in \mathcal{A}\}$ to 
\begin{align}
\underset{\boldsymbol{x}}{\text{maximize}} & \hspace{6pt} \sum_{k=1}^K p_{k0}(1-(1-\alpha_k)^{u_k}) \label{eq:original_obj_fun_1}
\end{align}
\begin{align}
\text{subject to} & \sum_{m : (m,k) \in \mathcal{A}} x_{mk} = u_k, \; \forall k \notag\\ 
& \sum_{k : (m,k) \in \mathcal{A}} x_{mk} = N_m, \; \forall m \notag\\
& x_{mk} \in \{0,1,\ldots,N_m\}, \; \forall (m,k) \in \mathcal{A} \notag
\end{align}
Note that the equality constraint for the agent search budgets can be included because of the monotonicity of the objective function.
The problem has a separable nonlinear objective function \eqref{eq:original_obj_fun_1} and linear  constraints.  Each of the individual functions in \eqref{eq:original_obj_fun_1} is concave when the variables $u_k$ are relaxed to be continuous. By introducing a few additional variables, we will be able to transform \eqref{eq:original_obj_fun_1} into linear form, as each of the individual functions in \eqref{eq:original_obj_fun_1} can be decomposed as 
\begin{align}
p_{k0}(1-(1-\alpha_k)^{u_k}) &= \underset{\{y_{kj}\}}{\text{maximize}} \hspace{6pt}  \sum_{j=1}^N p_{kj} y_{kj}\label{e:dac1} \\ 
\text{subject to} \hspace{6pt} & \sum_{j=1}^N y_{kj} = u_k ;\ 
 y_{kj} \in \{0,1\}, \; \forall j \notag
\end{align}

\section{A min-cost flow interpretation}
\label{sec:perspective}

We next map the above problem to a minimum cost network flow problem.  By using the transformation in \eqref{e:dac1}, we obtain the following integer programming problem:
\begin{align}
\underset{\boldsymbol{x}, \boldsymbol{y}}{\text{minimize}} \hspace{6pt} & \sum_{k=1}^K \sum_{j=1}^N - p_{kj} y_{kj} \notag\\
\text{subject to} \hspace{6pt} & \sum_{k : (m,k) \in \mathcal{A}} x_{mk} = N_m, \; m = 1,\ldots,M \label{eq:source_conserve}\\
& \sum_{m : (m,k) \in \mathcal{A}} x_{mk} = \sum_{j=1}^N y_{kj}, \; k=1,\ldots,K \label{eq:sink_conserve}\\
& \sum_{k=1}^K \sum_{j=1}^N y_{kj} = N; \ y_{kj} \in \{0,1\}, \; \forall k,j  \label{eq:global_conserve}\\
& x_{mk} \in \{0,1,\ldots,N_m\}, \; \forall (m,k) \in \mathcal{A} \notag
\end{align}

\begin{figure}[h]
\centering
\includegraphics[width=0.5\textwidth]{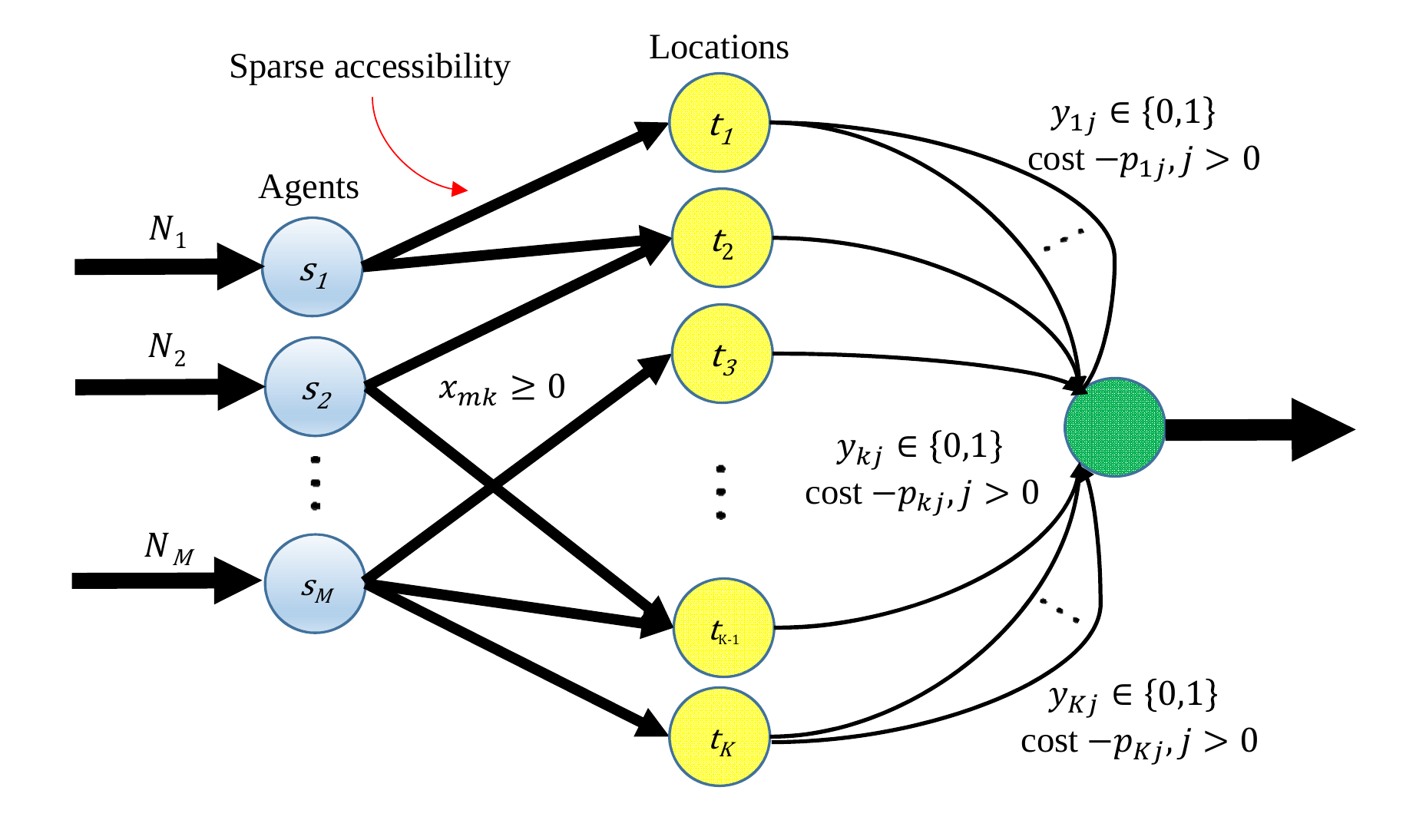}
\caption{Sparse multi-agent discrete search problem viewed as a min-cost network flow problem.}
\label{fig:flow_chart}
\end{figure}

A graphical representation of the constraints of this optimization problem is depicted in the directed multigraph in Figure \ref{fig:flow_chart}.  Source node $s_m$ with supply $N_m$ represents the $m$-th search agent.  Sink node $t_k$ represents the $k$-th location.  If agent $m$ can access location $k$, there is a directed arc from source node $s_m$ to sink node $t_k$.  

Let $x_{mk} \in \{0,1,\ldots,N_m\}$ denotes the flow from $s_m$ to $t_k$.  There is also a dummy global sink node.  The cost on each arc from source $s_m$ to sink $t_k$ is zero.  From $t_k$ to the global sink node, there are $N$ directed arcs, with cost $-p_{kj}$ on the $j$-th arc.  Denote the flow  on the $j$-th arc from $t_k$ to the global sink node by $y_{kj} \in \{0,1\}$.  

Equations \eqref{eq:source_conserve}, \eqref{eq:sink_conserve} and \eqref{eq:global_conserve} correspond to flow conservation at the source nodes $s_m$'s, the sink nodes $t_k$'s and the global sink node, respectively.   This integer linear programming problem is equivalent to the transformed problem presented in the previous section.  

From Fig. \ref{fig:flow_chart}, it is clear that the constraints in this integer program are unimodular, and the right-hand sides of the constraints are integers (the search efforts.)  Hence,  the optimal solutions of the linear programming relaxation are integer-valued.  Thus, we can obtain optimal solutions for the above problem using any network optimization package.  However, doing so requires explicit construction of the network in Fig. \ref{fig:flow_chart}, which requires creating a large number of additional arcs in the network. In our development below, we exploit the special structure of the network problem to develop a fast algorithm that avoids this extra construction. 

%
%
%

\subsection{Duality and complementary slackness } 

The Lagrangian of the primal is
\begin{small}
\begin{align*}
& L(\boldsymbol{x},\boldsymbol{y},\boldsymbol{d^s},\boldsymbol{d^t},\lambda) \\
= & \sum_{k=1}^K \sum_{j=1}^N - p_{kj} y_{kj} + \sum_{m=1}^M d^s_m (N_m - \sum_{k : (m,k) \in \mathcal{A}} x_{mk}) \\
& + \sum_{k=1}^K d^t_k (\sum_{m : (m,k) \in \mathcal{A}} x_{mk} - \sum_{j=1}^N y_{kj}) + \lambda (\sum_{k=1}^K \sum_{j=1}^N y_{kj} - N) 
\end{align*}
\end{small}

%

Here we let $d^s_m$, $d^t_k$, and $\lambda$ be the dual variables, a.k.a. the \textit{prices}, of source node $s_m$, sink node $t_k$ and the global sink node, respectively.  The dual problem is then:
\begin{small}
\begin{align}
\underset{\boldsymbol{d^s},\boldsymbol{d^t},\lambda}{\text{maximize}} \hspace{6pt} & \sum_{k=1}^K \sum_{j=1}^N \min\{0,  \lambda - d^t_k - p_{kj}\} + \sum_{m=1}^M d^s_m N_m - \lambda N \notag\\
\text{subject to} \hspace{6pt} & d^s_m \le d^t_k, \; \forall (m,k) \in \mathcal{A}  \label{eq:dual_fea_cond}
\end{align}
\end{small}

The complementary slackness conditions (see Theorem 9.4 of \cite{ahuja1993network}) are
\begin{align}
d^s_m < d^t_k & \Rightarrow x_{mk} = 0, \forall (m,k) \in \mathcal{A} \label{eq:comp_slack_1}\\
x_{mk} > 0 & \Rightarrow d^s_m = d^t_k, \forall (m,k) \in \mathcal{A} \label{eq:comp_slack_2}\\
d^t_k - \lambda < - p_{kj} & \Rightarrow y_{kj} = 0, \ \forall k, \forall j \label{eq:comp_slack_4}\\
0 < y_{kj} < 1 & \Rightarrow d^t_k - \lambda = -p_{kj}, \ \forall k, \forall j \label{eq:comp_slack_5}\\
d^t_k - \lambda > - p_{kj} & \Rightarrow y_{kj} = 1, \ \forall k, \forall j \label{eq:comp_slack_6}
\end{align}
By Theorem 9.4 of \cite{ahuja1993network}, if a primal feasible solution $\{\boldsymbol{x}, \boldsymbol{y}\}$ and a dual feasible solution $\{\boldsymbol{d^s}, \boldsymbol{d^s}, \boldsymbol{\lambda}\}$ satisfy the complementary slackness conditions \eqref{eq:comp_slack_1}-\eqref{eq:comp_slack_6}, then $\{\boldsymbol{x}, \boldsymbol{y}\}$ is also an optimal solution for the primal problem.  


\section{Fast algorithm for sparse multi-agent discrete search}
\label{sec:algorithm}

Here we propose a specialized algorithm for the problem of sparse multi-agent discrete search without false alarms.  Our algorithm is based on primal-dual techniques: it begins with a solution satisfying complementary slackness.  It continually increases the flow and modifies dual variables implicitly until primal feasibility is achieved, while maintaining complementary slackness at every iteration.  A notable feature of our algorithm is that dual variables can be tracked implicitly without computation, as will be shown below. 

















\vspace{10pt}
\noindent
\underline{\fontfamily{zi4}\selectfont Algorithm}
\vspace{5pt}

\noindent \textbf{Input:} $\mathcal{A}$, $N_m$, $p_{k0}$, $\alpha_k$, $\forall m, k$.

\noindent \textbf{Output:} $x_{mk}, \forall (m,k) \in \mathcal{A}$.

\noindent \textbf{Initialization:}
\begin{itemize}
\item[]
Compute and insert $-p_{k1}, k=1,\ldots,K$ into a min-oriented binary heap $H$ of size $K$.  Define the available supply at each source node $s_m$ as $R_m$.  Let $R_m = N_m$, $m=1,\ldots,M$.
\end{itemize}
\noindent
\textbf{Step 1:} Repeat until $H$ is empty:
\begin{itemize}
\item[] Extract the minimum $-p_{kj}$ from $H$:
\begin{itemize}
\item[a)]  If $t_k$ is marked as \textit{eliminated}, continue step 1.

\item[b)] Run {\fontfamily{zi4}\selectfont Assign-Extra-Demand($t_k$)} and let its returned value be $m$.

\item[c)] If $m>0$, decrement $R_m \leftarrow R_m-1$.  Let $y_{kj}=1$.  Compute and insert $-p_{k,j+1}$ into $H$. 

\item[d)] Else, identify $s_m$'s and $t_k$'s marked as \textit{visited} in {\fontfamily{zi4}\selectfont Assign-Extra-Demand} as an isolated group, and mark them as \textit{eliminated}.

\item[e)] Return to Step 1.
\end{itemize}
\end{itemize}
\noindent
\textbf{Step 2}: If $R_m = 0$ for all $m$, identify $s_m$'s and $t_k$'s that are never \textit{eliminated} as an isolated group. Terminate the algorithm. 

\vspace{10pt}
In the algorithm, we use subroutine {\fontfamily{zi4}\selectfont Assign-Extra-Demand} to search for an augmenting path from a sink node $t$ that requires assignment of an extra unit of demand to a source node that has available supply, and identifies isolated groups when it fails to discover such a path.  The subroutine is described as follows:  

\vspace{10pt}
\noindent
\underline{\fontfamily{zi4}\selectfont Assign-Extra-Demand($t$)}
\vspace{5pt}

\noindent
\textbf{Initialization:}  
\begin{itemize}
\item[] Begin with a queue $Q$ containing the single node $t$, and a set of \textit{visited} nodes $V$ initially empty.   Define available supply at each source $s_m$ as $R_m$.  Initialize an array of predecessors for all sources and sinks as $pred(s_m) = -1$, $pred(t_k) = -1$.  
\end{itemize}

\noindent
\textbf{Step 1}:  Remove top element in $Q$: 
\begin{itemize}
\item[a)] If element is a sink $t_k$, check each source $s_m$ such that $(m, k) \in \mathcal{A}$ and 
$s_m \not\in Q$ and $s_m \not\in V$ to $Q$.  If source $s_m$ has $R_m > 0$, an augmenting path has been found: set $t_k = pred(s_m)$ and go to \textbf{Step 3} below.  Else, add each such 
$s_m$ to $Q$ and set $pred(s_m) = t_k$.  Add $t_k$ to $V$ and continue.  

\item[b)] If element is a source $s_m$ and $R_m = 0$, find each sink $t_k$ such that $(m, k) \in \mathcal{A}$ and $x_{mk} > 0$ and $t_k \not\in Q$ and $t_k \not\in V$.  Add each such sink $t_k$ to $Q$, and set $pred(t_k) = s_m$.   Add $s_m$ to $V$ and continue.  

\end{itemize}
\noindent
\textbf{Step 2}:  If $Q$ is note empty, return to {\bf Step 1}.  Else, terminate the algorithm and mark all nodes in $V$ as \textit{eliminated}.  Shrink the graph.  Return -1.  

\noindent
\textbf{Step 3}: An augmenting path has been found. Set $m^*=m$.  Then, recursively, starting with $s = s_m$,  perform an augmentation as follows:  
\begin{itemize}
\item[i.] For source $s_m$, find $t_k = pred(s_m)$.  Increment $x_{mk} \leftarrow x_{mk}+1$.   

\item[ii.] For sink $t_k$, if $pred(t_k) = -1$, finish.  Else, let $s_m = pred(t_k)$ and decrement $x_{mk} \leftarrow x_{mk}-1$.  
\end{itemize}
\noindent
Terminate the subroutine  and return $m^*$.  
\vspace{10pt}  

If {\fontfamily{zi4}\selectfont Assign-Extra-Demand} fails to assign the extra unit of demand to any supply source, we identify all source and sink nodes marked as \textit{visited} during the search as an isolated group, and remove them from further consideration in the subsequent assignment problems by marking their status as \textit{eliminated}. 

Our algorithm has two notable features.  First, it exploits the implicit ordering of $p_{kj}$'s for the same location and computes $p_{kj}$'s on the fly only when they are needed.  This is assisted by maintaining a min-oriented heap.  Second, the idea of \textit{isolated groups} is essential for our algorithm.  At each iteration a $p_{kj}$ is popped out of the heap, creating one unit of demand at a sink.  The algorithm then tries to assign this extra unit of demand to a source with available supply.  When the subroutine {\fontfamily{zi4}\selectfont Assign-Extra-Demand} fails to find an augmenting path, we identify an isolated group, consisting of all the sink and source nodes marked as \textit{visited} during the breadth-first search, which no longer needs to be considered in the algorithm.    


Next we will show that the algorithm always terminates with finite time complexity, and analyze its time complexity.

\begin{theorem} \label{thm1}
The algorithm terminates in finite time, with complexity  is $O(N|\mathcal{A}|)$, where $N=\sum_{m=1}^M N_m$ is the total supply from all source nodes, and $|\mathcal{A}|$ is the cardinality of the set of all accessibility pairs.
\end{theorem}

\noindent The complete proof is included in the appendix.  Each iteration assigns an extra unit of search time from a source to a sink, or eliminates one or more sinks from future searches. The worst-case complexity of finding a source with available supply is $O(|\mathcal{A}|)$, which forms the basis for the complexity bound.  Note that the complexity of a fast min-cost flow algorithm is  $O\big((|\mathcal{A}|+N)(|\mathcal{A}|+N+(M+K)\log(M+K))\log(\max_m N_m)\big)$ \cite{ahuja1993network}, which is significantly larger than the complexity above.   
Our algorithm also works for the problem considered in \cite{song2004discrete}, in which every agent can access every location.  Assuming that there are $M$ agents, the amortized complexity of finding an augmenting path is $O(1)$, so the complexity of our algorithm becomes $O(N \log(K))$.     In contrast, the algorithm of \cite{song2004discrete} requires $O(K N/M \log(KN/M))$, with the assumption that $K > M$, slower than our algorithm. 


\noindent
\begin{theorem} \label{thm2}
The algorithm constructs an optimal multi-agent search allocation. 
\end{theorem}

The proof in the appendix shows that there exists a set of dual prices that satisfy complementary slackness along with the primal feasible allocation constructed by the algorithm, thus establishing optimality. 

\section{Probabilies of Detection Depend on Both Search Locations and Search Agents}
\label{sec:discussion}

Now we consider the case where the probabilities of detection depend on both search agents and locations.  Denote the detection probability at location $k$ using agent $m$ as $\alpha_{mk}$.  Suppose that we allocate agent $m$ to search location $k$ for $x_{mk}$ times.  Let $\boldsymbol{x}_k=(x_{1k},\ldots,x_{Mk})$, so the probability of finding the object at location $k$ under allocation $\boldsymbol{x}_k, k=1,\ldots,K$ is $P_k(\boldsymbol{x}_k) = p_{k0}(1-\prod_{m=1}^M (1-\alpha_{mk}) ^ {x_{mk}})$.  Then the problem is to find $\boldsymbol{x}=\{\boldsymbol{x}_k\}$ to 

\begin{align}
\underset{\boldsymbol{x}}{\text{maximize}} & \hspace{6pt} \sum_{k=1}^K  p_{k0}\left(1-\prod_{m=1}^M (1-\alpha_{mk}) ^ {x_{mk}}\right) \label{eq:np_hard_obj_fun_1}\\
\text{subject to} 
& \sum_{k : (m,k) \in \mathcal{A}} x_{mk} \le N_m, \; \forall m \notag\\
& x_{mk} \in \{0,1,\ldots,N_m\}, \; \forall (m,k) \in \mathcal{A} \notag
\end{align}

This problem is known to be NP-hard \cite{lloyd1986weapons, ahuja2007wta}.  Approximation schemes such as branch-and-bound and local search have been proposed recently for these problems \cite{ahuja2007wta}.  

An alternative formulation of the problem is based on monotone submodular objective functions and matroid constraints.  Without loss of generality, assume that we separate the search effort of each of the $M$ agents by giving them unique identities, so that we have an equivalent problem with $N$ total agents, and the constraint that $N_{n} = 1$ for each agent $n$.  Let $\mathcal{E}$ be the set of all accessibility pairs after separation, and let $\alpha_{nk}$ be the corresponding detection probability at location $k$ using agent $n$. 

Denote by the subset $S_n$ of $\mathcal{E}$ as $S_n = \{(n,k) | k \in \{1, \ldots, K\}, (n,k) \in \mathcal{E} \}$, the set of accessibility pairs of agent $n$.  Note that $S_n \cap S_{n'} = \emptyset$ if $n \ne n'$, and 
$$ \mathcal{E} = \bigcup_{n = 1}^N S_n$$
so that $S_1,  \ldots, S_N$ is a partition of $\mathcal{E}$.  A feasible assignment $S$ from the agents to the locations is a subset $S \subset \mathcal{E}$ such that 
$$ |S \cap S_n| \le 1, \quad n = 1, \ldots, N$$  

The value of this feasible assignment is given by  
\begin{equation} 
\label{e:dac4}
f(S) = \sum_{k=1}^K  p_{k0}\left(1-\prod_{(n,k) \in S}(1-\alpha_{nk}) \right)
 \end{equation}

$f(S)$ has the following properties:

\begin{lemma} \label{lemma_fs}
$f(S)$ is increasing and submodular in $S$.
\end{lemma}

The submodularity of $f(S)$ is proved by showing that it has the diminishing return property \cite{fisher_nem_submod_2}.

Let the collection of feasible assignments be $\mathcal{I}$.  Then $(\mathcal{E},\mathcal{I})$ has the following property:

\begin{lemma} \label{lemma_matroid}
The pair $(\mathcal{E},\mathcal{I})$ is a matroid. 
\end{lemma}
This can be proved by showing that $(\mathcal{E},\mathcal{I})$ satisfies both the hereditary and the augmentation properties \cite{welsh2010matroid}.  The proofs of Lemma \ref{lemma_fs}  and \ref{lemma_matroid} are included in the appendix.

Thus, the equivalent problem of maximizing \eqref{e:dac4} subject to $S \in \mathcal{I}$ is a submodular maximization problem subject to a matroid constraint.  It is known that the greedy algorithm produces a guarantee of 1/2-approximation \cite{fisher_nem_submod_2}.  This approximation ratio can be improved to $(1-1/e)$ by using a randomized algorithm with pipage rounding and a continuous greedy process \cite{calinescu2011}.  

We briefly describe a simple greedy algorithm here. This algorithm works directly with the original version of the problem and does not require separation of the sources into individual nodes, thus working with the much smaller network representation $\mathcal{A}$ instead of $\mathcal{E}$.
  Initialize $\boldsymbol{x}$ to be all zero.  Define available supply at each source as $R_m$ and let $R_m=N_m$.  Insert the triples $(\max_{\{m:(m,k) \in \mathcal{A}\}} \alpha_{mk} p_{k0}, k, \text{arg}\max_{\{m:(m,k) \in \mathcal{A}\}} \alpha_{mk}), \forall k$ into a max heap $\mathcal{H}$ ordered by the first value of the triples.  While there are elements in the heap, perform the following steps:
\begin{itemize}
\item Remove top element $(V, k, m')$ in the heap. 
\item If $R_{m'} > 0$,
\begin{itemize}
\item  Increment $x_{m'k} = x_{m'k}+1$ and decrement $R_{m'} = R_{m'} - 1$.
\item Define $V' = V(1-\alpha_{m'k})$ and insert $(V',k,m')$ into heap $\mathcal{H}$.
\end{itemize}
\item Else $R_{m'} = 0$: 
\begin{itemize}
\item Find $m'' = \argmax_{\{m : R_m > 0, (m,k) \in \mathcal{A}\}} \alpha_{mk}$.  
\item If $m''$ exists, then compute $V' = V \frac{\alpha_{m''k}}{\alpha_{m'k}}$ and insert $(V', k, m'')$ into heap $\mathcal{H}$.  If no $m''$ exists, continue to the next element in the heap.
\end{itemize}
\end{itemize}
This algorithm terminates with a greedy allocation once all the supplies $R_m$ decrease to 0.

\section{Experiments}
\label{sec:experiments}

In this section, we perform two experiments to validate the time complexity of our algorithm as derived in Theorem \ref{thm1}, and compare its running time with a classic min-cost flow algorithm, the capacity scaling algorithm.  

\begin{figure}[thpb]
\centering
\includegraphics[width=0.3\textwidth]{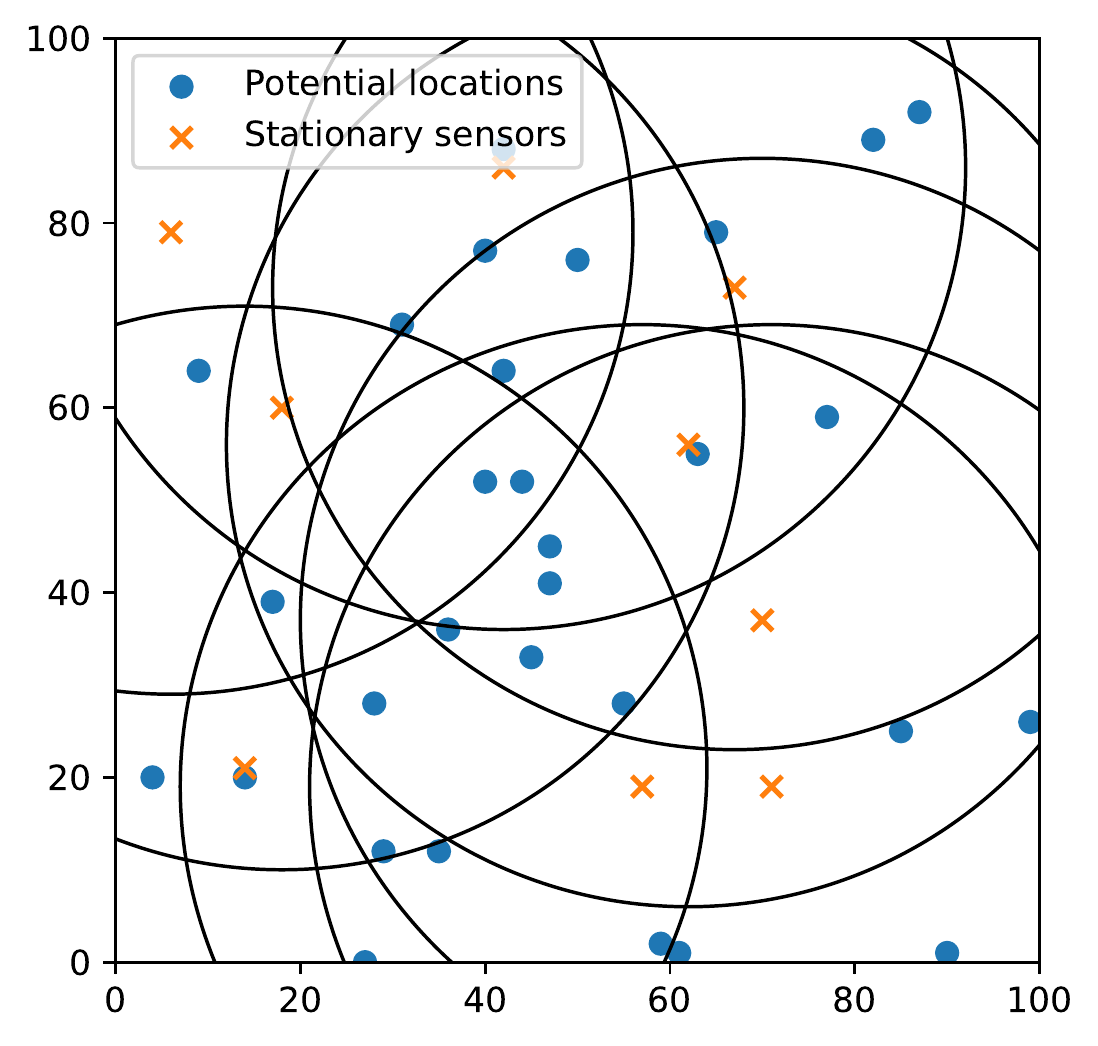}
\caption{An example spatial field with 9 stationary sensors and 30 potential locations.  The circle around each sensor represents its search coverage range.}
\label{fig:spatial_field_example}
\end{figure}

We conduct the experiments using Python 3 on a laptop computer with Intel i7-4600M processor and 8GB RAM.  We consider the physical setting where multiple stationary sensors with limited sensing range are used to search multiple potential locations in a 2-D spatial field.  We randomly generate the potential locations and the deployed sensors while making sure that each potential location lies within the search radius of at least one sensor and each sensor can search at least one potential location.  An example spatial field is in Figure \ref{fig:spatial_field_example}.  We randomly generate the prior probability $p_{k0}$ and the probability of detection $\alpha_k$ for each potential location $k$. Then we build a graph as described in Section \ref{sec:perspective}.  The graph data structure we use is provided by NetworkX \cite{networkxPaper}, a popular Python graph library.  The capacity scaling algorithm that we compare our algorithm with is also provided by NetworkX.

\begin{figure}[thpb]
\centering
\includegraphics[width=0.4\textwidth]{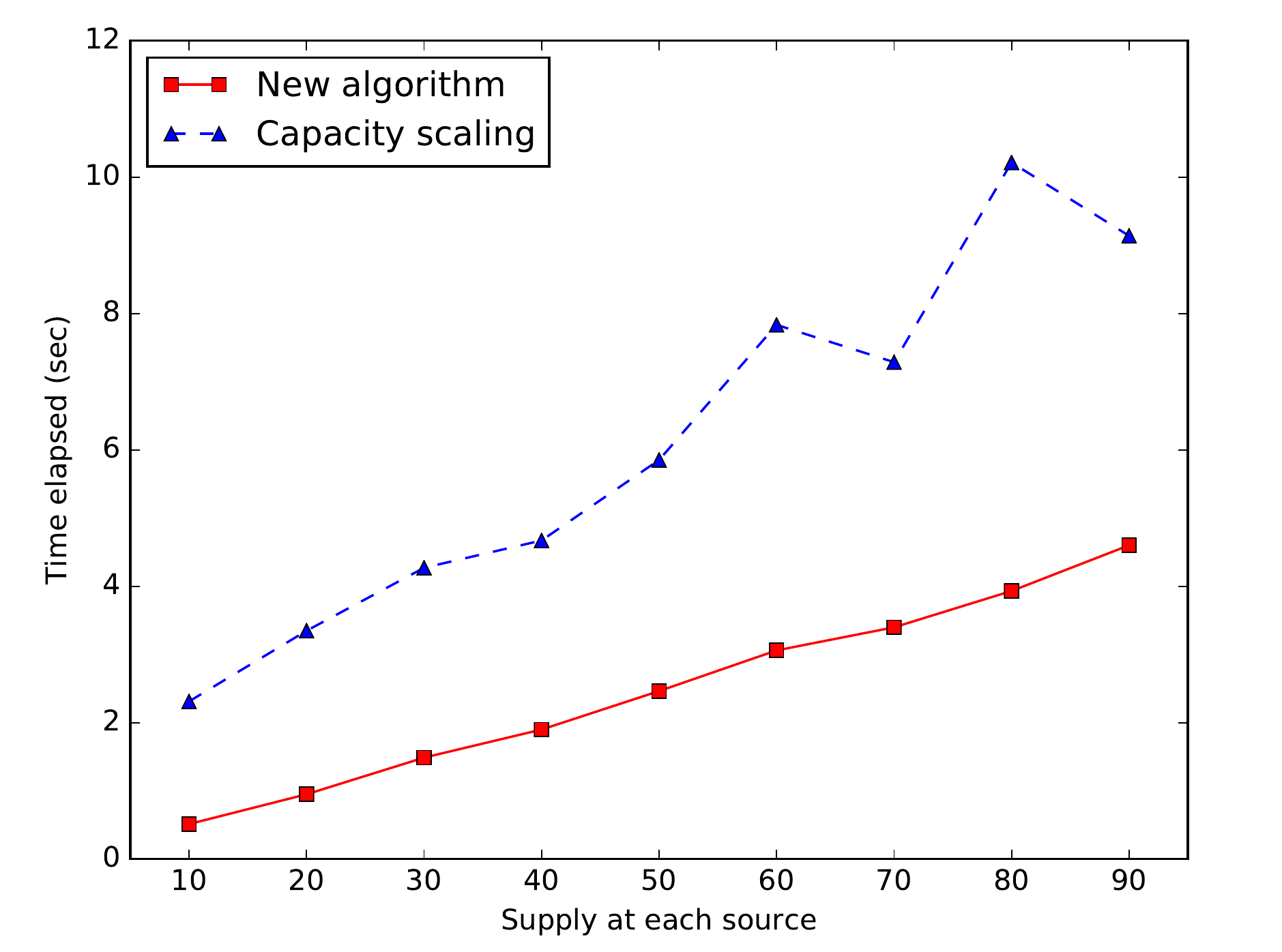}
\caption{Comparison of running time between the new algorithm and the capacity scaling algorithm, for different supply $N_{each}$ at each source node. $\#\text{sources}=100$, $\#\text{sinks}=1000$.}
\label{fig:supply_UAV}
\end{figure}

In the first experiment, we study the effect of source supply (i.e. sensor budget) on the running time of the two algorithms.  For simplicity, we let each sensor have the same budget $N_{each}$.  We randomly generate a spatial field, with 100 sensors and 1000 potential locations.  We let the search radius of each sensor be 15.  Then we vary the budget $N_{each}$ of each sensor from 10 to 90 with a step of 10, and compare the running time of the two algorithms.  The result is plotted in Figure \ref{fig:supply_UAV}. 

\begin{figure}[thpb]
\centering
\includegraphics[width=0.4\textwidth]{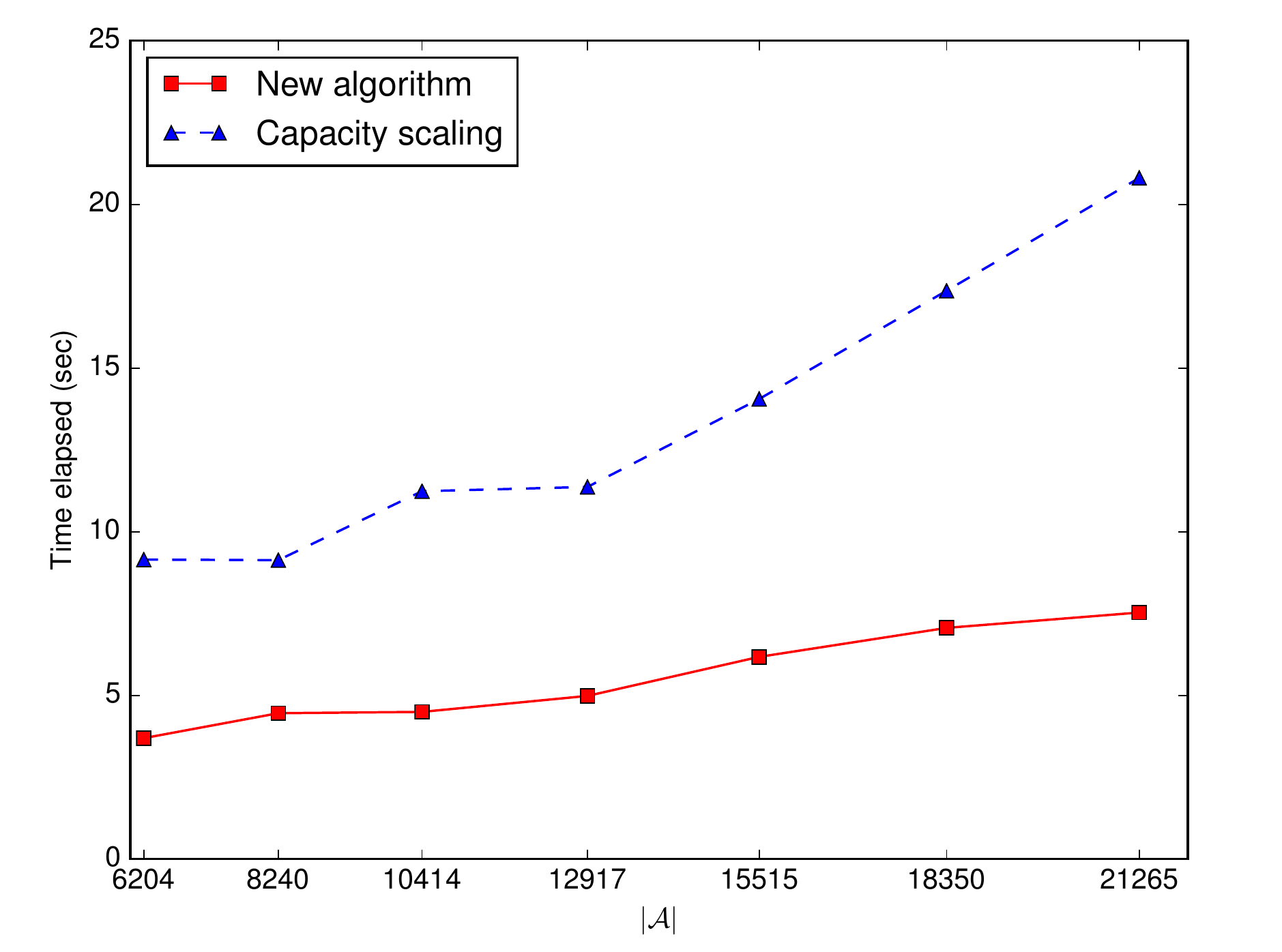}
\caption{Comparison of running time between the new algorithm and the capacity scaling algorithm, for different $|\mathcal{A}|$. $\#\text{sources}=100$, $\#\text{sinks}=1000$, $N_{each}=50$.}
\label{fig:sparsity_UAV}
\end{figure}

In the second experiment, we study the effect of graph sparsity, or the number of arcs between sources and sinks $|\mathcal{A}|$, on the running time of the two algorithms.  We randomly generate a spatial field, with 100 sensors and 1000 potential locations.  We let each sensor have a budget $N_{each}=50$.  Then we vary the search radius of each sensor from 15 to 30 with a step of 2.5, and count the corresponding $|\mathcal{A}|$ for different search radiuses.  We compare the two algorithms and plot their running time with respect to $|\mathcal{A}|$ in Figure \ref{fig:sparsity_UAV}. 

From the figures, we can see that our specialized algorithm scales linearly with respect to $N_{each}$ and $|\mathcal{A}|$.  These experimental results have validated the time complexity analysis in Theorem \ref{thm1}.  Besides, for both experiments, our algorithm runs multiple times faster than the capacity scaling algorithm.  This demonstrates that our specialized algorithm is highly efficient for solving the problem of sparse multi-agent search without false alarms studied in this paper.  

\section{Conclusion}
\label{sec:conclusion}

In this paper, we investigated the problem of sparse multi-agent discrete search without false alarms.  We provided a novel perspective of viewing the problem as a min-cost flow problem, and gave the optimality conditions for the agent schedule to be optimal.  We proposed a fast specialized algorithm for solving the problem.  The algorithm uses a heap to create one unit of demand at a sink node at each iteration, then finds an augmenting path to assign the demand.  When it fails to do so, we identify and remove an isolated group of nodes from further consideration.  We proved that the algorithm can always find an optimal solution in finite time, and analyzed the time complexity of the algorithm.  We performed experiments to compare our specialized algorithm with a general min-cost flow algorithm, the capacity-scaling algorithm.  

We also addressed the problem where detection performance depends on both location and agent, which is known to be NP-Hard.  We showed that the problem can be reduced to a submodular maxization problem over a matroid, and provided an approximate algorithm with guaranteed performance. 

There are several directions in which the paper can be extended.  One direction is to study multi-agent whereabouts search without false alarms where the objective is to find an optimal agent schedule to maximize the probability of correctly stating where the object is.  Another direction is to consider when one search action of different agents cost different amounts of effort.  Optimal strategies may not exist for this problem; techniques such as branch and bound may be needed.  

\section*{Appendix}
\label{sec:appendix}

{\bf Proof of Theorem \ref{thm1}:}

We first show finite time termination.  Each iteration of the algorithm either successfully assigns one unit of demand to some supply source, or identifies and eliminates an isolated group which consists of at least one source node and one sink node.  Given the connectivity assumptions, there is a feasible assignment of all source supply to sink nodes.  Since we have a finite number of source and sink nodes, and a finite integer number of supply, the algorithm will eventually terminate.  Upon the termination of the algorithm, due to the assumption that each location is accessible by at least one agent and each agent can access at least one location, there will be no remaining unassigned supply. 

In terms of computation, the overhead for constructing the binary heap is $O(K)$.  The worst-case running time for heap insertion and deletion is $O(\log(K))$.  The worst-case running time for finding an augmenting path is $O(|\mathcal{A}|)$ since we might have to explore every arc between the source and sink nodes.  At each iteration of the algorithm, the algorithm may successfully assign one unit of supply and do heap insertion and deletion, or fail to find an augmenting path.  There are $N=\sum_{m=1}^M N_m$ units of supply to assign, and there are at most $\min(M,K)$ failures since each failure eliminates at least one source and one sink.  Thus, the overall running time is $O\big(K + N (\log(K)+|\mathcal{A}|) + \min(M,K) |\mathcal{A}|\big) = O\big(N|\mathcal{A}|\big)$ since $\max(M,K) \le |\mathcal{A}|$ and $\min(M,K) \le M \le N$. 
$\QED$

\vspace{10pt}

{\bf Proof of Theorem \ref{thm2}:}

The algorithm constructs a primal feasible solution $\{\boldsymbol{x}, \boldsymbol{y}\}$ which satisfies all the constraints of the primal problem.  As discussed before, to show that $\{\boldsymbol{x}, \boldsymbol{y}\}$ is also optimal, we only need to show that there exist prices that satisfy the dual feasibility condition \eqref{eq:dual_fea_cond} and the complementary slackness conditions \eqref{eq:comp_slack_1}-\eqref{eq:comp_slack_6}.  

We will assign prices as follows:  For each isolated group of source and sink nodes, assume that $-p_{k,j^*(k)}$ is the value popped from the min-oriented binary heap in the algorithm that led to the failed augmenting path search (or a successful augmenting path search and the termination of the algorithm since all supplies are used up), and resulted in the formation of the isolated group.  Set the prices of source nodes $\{s_m\}$ and sink nodes $\{t_k\}$ in this isolated group to $d^s_m = d^t_k = -p_{k,j^*}$.  

When the algorithm terminates, there will be remaining source and sink nodes that are never marked as \textit{eliminated}. Set the prices of these source nodes $\{s_m\}$ and sink nodes $\{t_k\}$ to $d^s_m = d^t_k = -p_{k', j^*(k')}$, the last value extracted from the min-oriented binary heap before the algorithm terminates.  

Additionally, we let the price $\lambda$ of the global sink node be 0. 

Note that the price $d^t_k$ of sink node $t_k$ is equal to the last (and largest) extracted value of the isolated group to which $t_k$ belongs.  So for the $j$-th arc from sink node $t_k$ to the global sink node, if $-p_{kj} < d^t_k$, it implies that $-p_{kj}$ has already been extracted from the min-oriented heap and the flow on the arc $y_{kj}$ has been assigned to 1.  Since $\lambda=0$, we have $d^t_k - \lambda > - p_{kj} \Rightarrow y_{kj} = 1$; \eqref{eq:comp_slack_6} is satisfied.  Similarly, if $-p_{kj} > d^t_k$, it implies that $-p_{kj}$ has not been extracted before the isolated group of $t_k$ is eliminated, which is equivalent to $y_{kj}=0$.  Thus, $d^t_k - \lambda < - p_{kj}  \Rightarrow y_{kj} = 0$; \eqref{eq:comp_slack_4} is satisfied.  In addition, since every $y_{kj}$ must be either 0 or 1, \eqref{eq:comp_slack_5} is automatically satisfied.  

For the arc from source node $s_m$ to sink node $t_k$,  if $d^s_m < d^t_k$, it implies that the isolated group of $s_m$ is identified and eliminated before the isolated group of $t_k$.  So $x_{mk}=0$, since otherwise {\fontfamily{zi4}\selectfont Assign-Extra-Flow} would visit $t_k$ and eventually put it in the same isolated group with $s_m$, leading to $d^s_m = d^t_k$.   Once $s_m$ is eliminated from the graph, so is $x_{mk}$, thus $x_{mk}$ will remain 0 and never be updated.  Thus, $d^s_m < d^t_k \Rightarrow x_{mk} = 0$; \eqref{eq:comp_slack_1} is satisfied.  On the other hand, if $x_{mk} > 0$, {\fontfamily{zi4}\selectfont Assign-Extra-Flow} will always visit $t_k$ after visiting $s_m$, or visit $s_m$ after visiting $t_k$, and put them into the same isolated group when it fails, leading to $d^s_m=d^t_k$.  Thus, $x_{mk} > 0  \Rightarrow d^s_m = d^t_k$; \eqref{eq:comp_slack_2} is satisfied.  Finally, there will not exist $d^s_m$ and $d^t_k$ such that $d^s_m > d^t_k$. This is because $d^s_m > d^t_k$ implies that the isolated group of $t_k$ is identified and eliminated before the isolated group of $s_m$.  However, this cannot happen since if {\fontfamily{zi4}\selectfont Assign-Extra-Flow} visits $t_k$, it will always visit $s_m$, no matter what value $x_{mk}$ is, as long as $(m,k) \in \mathcal{A}$.  Thus, the dual feasibility condition \eqref{eq:dual_fea_cond} is satisfied.  The proof is now complete. $\QED$

\vspace{10pt}

{\bf Proof of Lemma \ref{lemma_fs}:}

Given two feasible assignments $S' \subset S$, it is easy to see that $f(S') \le f(S)$, because each term of the sum in \eqref{e:dac4} is smaller in $f(S')$, as all the products in $S'$ are in $S$ and $S$ has additional products.  Thus, $f(S)$ is increasing.  

Furthermore, $f(S)$ is a submodular function because given $S' \subset S$ and $e = (\tilde{m}, \tilde{k}) \notin S$
\begin{align*}
f(S' \cup \{e\}) - f(S') &=  p_{k0} \alpha_{\tilde{m},\tilde{k}} \prod_{(m,k) \in S'}(1-\alpha_{mk})  \\
&\ge   p_{k0} \alpha_{\tilde{m},\tilde{k}} \prod_{(m,k) \in S}(1-\alpha_{mk}) \\
& =  f(S \cup \{e\}) - f(S) \hspace{60pt}\QED
\end{align*}

\vspace{10pt}

{\bf Proof of Lemma \ref{lemma_matroid}:}

For any feasible assignment $S \subset \mathcal{I}$, any $S' \subset S$ is also a feasible assignment, satisfying the hereditary property of independent sets required for a matroid.  

Furthermore, assume we have two feasible assignments $S, S'$ such that $|S'| < |S|$.  Then, there exists a source $n$ such that $(n,k) \in S$, and $(n,k') \notin S'$ for all $k' = 1, \ldots, K$.  This implies that $S' \cup \{(n,k)\}$ is also a feasible assignment, satisfying the augmentation property and establishing that $(\mathcal{E},\mathcal{I})$ is a matroid. $\QED$

\bibliographystyle{IEEEtran}
\bibliography{references}

\begin{thebibliography}{10}
\providecommand{\url}[1]{#1}
\csname url@samestyle\endcsname
\providecommand{\newblock}{\relax}
\providecommand{\bibinfo}[2]{#2}
\providecommand{\BIBentrySTDinterwordspacing}{\spaceskip=0pt\relax}
\providecommand{\BIBentryALTinterwordstretchfactor}{4}
\providecommand{\BIBentryALTinterwordspacing}{\spaceskip=\fontdimen2\font plus
\BIBentryALTinterwordstretchfactor\fontdimen3\font minus
  \fontdimen4\font\relax}
\providecommand{\BIBforeignlanguage}[2]{{%
\expandafter\ifx\csname l@#1\endcsname\relax
\typeout{** WARNING: IEEEtran.bst: No hyphenation pattern has been}%
\typeout{** loaded for the language `#1'. Using the pattern for}%
\typeout{** the default language instead.}%
\else
\language=\csname l@#1\endcsname
\fi
#2}}
\providecommand{\BIBdecl}{\relax}
\BIBdecl

\bibitem{Koopman1946}
B.~O. Koopman, \emph{Search and Screening}.\hskip 1em plus 0.5em minus
  0.4em\relax Operations Evaluation Group, Office of the Chief of Naval
  Operations, Navy Department, 1946.

\bibitem{StoneBook}
L.~D. Stone, \emph{Theory of Optimal Search}.\hskip 1em plus 0.5em minus
  0.4em\relax Academic Press, 1975.

\bibitem{ahlswede1987search}
R.~Ahlswede and I.~Wegener, \emph{Search Problems}.\hskip 1em plus 0.5em minus
  0.4em\relax Wiley, 1987.

\bibitem{castanon1995optimal}
D.~A. Casta{\~n}{\'o}n, ``Optimal search strategies in dynamic hypothesis
  testing,'' \emph{IEEE Trans. on Systems, Man, and Cybernetics}, vol.~25,
  no.~7, pp. 1130--1138, 1995.

\bibitem{song2004discrete}
N.-O. Song and D.~Teneketzis, ``Discrete search with multiple sensors,''
  \emph{Math. Methods of Oper. Res.}, vol.~60, no.~1, pp. 1--13, 2004.

\bibitem{ding2015optimal}
H.~Ding and D.~A. Casta{\~n}{\'o}n, ``Optimal solutions for classes of adaptive
  search problems,'' in \emph{Proc. IEEE Conf. Decision and Control}, Osaka,
  Japan, 2015.

\bibitem{ding2015multi}
------, ``Multi-object two-agent coordinated search,'' in \emph{Proc. Int'l
  Conf. Complex Systems Engineering}, Storrs, CT, 2015.

\bibitem{pattipati1990application}
K.~R. Pattipati and M.~G. Alexandridis, ``Application of heuristic search and
  information theory to sequential fault diagnosis,'' \emph{IEEE Trans. on
  Systems, Man, and Cybernetics}, vol.~20, no.~4, pp. 872--887, 1990.

\bibitem{acc16}
H.~Ding, E.~Cristofalo, J.~Wang, D.~A. Casta{\~n}{\'o}n, E.~Montijano,
  V.~Saligrama, and M.~Schwager, ``A multi-resolution approach for discovery
  and 3{D} modeling of archaeological sites using satellite imagery and a
  {UAV}-borne camera,'' in \emph{Proc. American Control Conf.}, Boston, MA,
  2016.

\bibitem{ding2016fast}
H.~Ding and D.~A. Casta{\~n}{\'o}n, ``Fast algorithms for {UAV} tasking and
  routing,'' in \emph{Proc. IEEE Conf. Control Applications}, 2016.

\bibitem{Benkoski1991}
S.~J. Benkoski, M.~G. Monticino, and J.~R. Weisinger, ``A survey of the search
  theory literature,'' \emph{Naval Research Logistics}, vol.~38, no.~4, pp.
  469--494, 1991.

\bibitem{Davidsbook}
A.~O. Hero, D.~A. Casta{\~n\'o}n, D.~Cochran, and K.~Kastella,
  \emph{Foundations and Applications of Sensor Management}.\hskip 1em plus
  0.5em minus 0.4em\relax Springer, 2008.

\bibitem{markovObj}
D.~C. Hitchings and D.~A. Casta{\~n}{\'o}n, ``Sensor control for search and
  identification of markov objects,'' \emph{Proc. IEEE Conf. Decision and
  Control and European Control Conf.}, December 2011.

\bibitem{tognetti1968}
K.~P. Tognetti, ``An optimal strategy for a whereabouts search,''
  \emph{Operations Research}, vol.~16, no.~1, pp. 209--211, 1968.

\bibitem{kadane1971whereabouts}
J.~B. Kadane, ``Optimal whereabouts search,'' \emph{Operations Research},
  vol.~19, no.~4, pp. 894--904, 1971.

\bibitem{wtaReport}
D.~A. Casta{\~n}{\'o}n, ``Advanced weapon-target assignment algorithm,''
  ALPHATECH, Inc., Burlington, MA, Tech. Rep., 1987.

\bibitem{lloyd1986weapons}
S.~P. Lloyd and H.~S. Witsenhausen, ``Weapons allocation is np-complete,'' in
  \emph{Proc. Summer Computer Simulation Conference}, 1986.

\bibitem{ahuja2007wta}
R.~K. Ahuja, A.~Kumar, K.~C. Jha, and J.~B. Orlin, ``Exact and heuristic
  algorithms for the weapon-target assignment problem,'' \emph{Operations
  Research}, vol.~55, no.~6, pp. 1136--1146, 2007.

\bibitem{ahuja1993network}
R.~K. Ahuja, T.~L. Magnanti, and J.~B. Orlin, \emph{Network Flows: Theory,
  Algorithms, and Applications}.\hskip 1em plus 0.5em minus 0.4em\relax
  Prentice hall, 1993.

\bibitem{fisher_nem_submod_2}
M.~L. Fisher, G.~L. Nemhauser, and L.~A. Wolsey, ``An analysis of
  approximations for maximizing submodular set functions —- {II},''
  \emph{Math. Prog. Study}, vol.~8, pp. 73--87, 1978.

\bibitem{welsh2010matroid}
D.~J.~A. Welsh, \emph{Matroid Theory}.\hskip 1em plus 0.5em minus 0.4em\relax
  Academic Press, 1976.

\bibitem{calinescu2011}
G.~Calinescu, C.~Chekuri, M.~P{\'a}l, and J.~Vondr{\'a}k, ``Maximizing a
  monotone submodular function subject to a matroid constraint,'' \emph{SIAM J.
  Computing}, vol.~40, no.~6, pp. 1740--1766, 2011.

\bibitem{networkxPaper}
D.~A. Schult and P.~Swart, ``Exploring network structure, dynamics, and
  function using {N}etwork{X},'' in \emph{Proc. Python Sci. Conf.}, 2008.

\end{thebibliography}
\end{document}